\documentclass[11pt,reqno]{amsart}
\usepackage{indentfirst, amssymb,amsmath,mathrsfs,amsthm,enumerate}
\usepackage[bookmarksnumbered, colorlinks, plainpages]{hyperref}
\newtheorem*{theo2A}{Theorem 2.A}
\newtheorem*{theo2B}{Theorem 2.B}
\newtheorem*{theo2C}{Theorem 2.C}
\newtheorem*{theo2D}{Theorem 2.D}

\newtheorem*{theo3A}{Theorem 3.A}
\newtheorem*{theo3B}{Theorem 3.B}
\newtheorem*{theo3C}{Theorem 3.C}
\newtheorem*{theo3D}{Theorem 3.D}

\newtheorem{ques}{Question}[section]

\newtheorem*{cor A}{Corollary A}
\newtheorem*{cor B}{Corollary B}

\newtheorem{theo}{Theorem}[section]
\newtheorem{lem}{Lemma}[section]
\newtheorem{cor}{Corollary}[section]

\newtheorem{rem}{Remark}[section]

\newcommand{\ol}{\overline}
\newcommand{\be}{\begin{equation}}
\newcommand{\ee}{\end{equation}}
\newcommand{\beas}{\begin{eqnarray*}}
\newcommand{\eeas}{\end{eqnarray*}}
\newcommand{\bea}{\begin{eqnarray}}
\newcommand{\eea}{\end{eqnarray}}

\numberwithin{equation}{section}
\begin{document}
\title[U\MakeLowercase{niqueness of First Derivatives}.....]{\LARGE U\huge\MakeLowercase{niqueness of First Derivatives and Differences in Meromorphic Functions and the Characterization of Entire Function Periodicity}}
\date{}
\author[A. B\MakeLowercase{anerjee}, S. M\MakeLowercase{ajumder} \MakeLowercase{and} N. S\MakeLowercase{arkar}]{A\MakeLowercase{bhijit} B\MakeLowercase{anerjee}, S\MakeLowercase{ujoy} M\MakeLowercase{ajumder}$^\ast$ \MakeLowercase{and} N\MakeLowercase{abadwip} S\MakeLowercase{arkar}}
\address{ Department of Mathematics, University of Kalyani, West Bengal 741235, India.}
\email{abanerjee\_kal@yahoo.co.in}
\address{Department of Mathematics, Raiganj University, Raiganj, West Bengal-733134, India.}
\email{sm05math@gmail.com, sjm@raiganjuniversity.ac.in}
\address{Department of Mathematics, Raiganj University, Raiganj, West Bengal-733134, India.}
\email{naba.iitbmath@gmail.com}

\renewcommand{\thefootnote}{}
\footnote{*\emph{Corresponding Author}: Sujoy Majumder.}

\renewcommand{\thefootnote}{\arabic{footnote}}
\setcounter{footnote}{0}

\begin{abstract}  The objective of the paper is twofold. The first objective is to study the uniqueness problem of meromorphic function $f(z)$ when $f^{(1)}(z)$ shares two distinct finite values $a_1$, $a_2$ and $\infty$ CM with $\Delta_cf(z)$. In this context, we provide a result that resolves the open problem posed by Qi et al. [Comput. Methods Funct. Theory, 18 (2018), 567-582] for the case when hyper order of the function is less than $\infty$. The second objective is to establish sufficient conditions for the periodicity of transcendental entire functions. In this direction, we obtain a result that affirms the question raised by Wei et al. [Anal. Math., 47 (2021), 695-708.]

\end{abstract}
\thanks{Typeset by \AmS -\LaTeX}
\maketitle

\section{{\bf Introduction}}
In recent years, the Nevanlinna value distribution of difference operators, as developed by Halburd and Korhonen \cite{3a} and Chiang and Feng \cite{1}, has played a crucial role in studying meromorphic solutions to complex difference equations and the uniqueness of complex differences. Numerous authors (see \cite{BM1}, \cite{BM2}, \cite{BM3}, \cite{DM1}, \cite{FLSY}, \cite{6}, \cite{7}, \cite{HF1}, \cite{M1}, \cite{MD},  \cite{MD1}, \cite{MP1}, \cite{MS1}, \cite{MS2}, \cite{MSP1}, \cite{10}, \cite{9a}) are currently investigating the value-sharing issues of meromorphic functions in connection to their difference operators, demonstrating the significant interest in this field.

\smallskip
We assume that the reader is familiar with standard notation and main results of Nevanlinna Theory (see \cite{14}). We shall call a meromorphic function $a(z)$ a small function of $f(z)$ if $T(r,a)=S(r,f)$, where $S(r, f)$ denotes any quantity satisfying the condition $S(r, f) = o(T(r, f))$ as $r\to \infty$ possibly outside of an exceptional set of finite logarithmic measure. Also we denote by $S(f)$ the set of functions which are small compared to $f$. In the paper $\rho(f)$ and $\rho_1(f)$ stand for the order and the hyper-order of a meromorphic function $f$ respectively. As usual, the abbreviation CM means counting multiplicities.

\smallskip 
We define an $\varepsilon$-set to be a countable union of discs (see \cite{1a})
\beas E=\sideset{}{_{j=1}^{\infty}}{\bigcup} B(b_j,r_j)\;\;\text{such}\;\text{that}\;\lim\limits_{j\to\infty}|b_j|=\infty\;\;\text{and}\;\;\sideset{}{_{j=1}^{\infty}}{\sum}r_j/|b_j|<\infty. \eeas

Here $B(a, r)$ denotes the open disc of center $a$ and radius $r$, and $S(a, r)$ will denote the corresponding boundary circle. Note that if $E$ is an $\varepsilon$-set, then the set of $r\geq 1$ for which the circle $S(0, r)$ meets $E$ has a finite logarithmic measure.

\smallskip
The paper is organized as follows. In Section 2, we study the uniqueness problem of meromorphic function $f$ when $f^{(1)}$ and $\Delta_cf$ share three values. In this section, we obtain a result which solves the open problem of Qi et al. \cite{10} for the case when $\rho_1(f)<\infty$. In Section 3, we study the sufficient condition for the periodicity of transcendental entire function and obtain a result which gives an affirmative answer of the question of Wei et al. \cite{11b}.

\section{\bf{Uniqueness problem of $f^{(1)}$ and $\Delta_cf$ sharing three values}}

In 1977, Rubel and Yang \cite{11} were the first to consider the uniqueness of an entire function and its derivative. They proved the following result:

\begin{theo2A}\cite{11} Let $f$ be a non-constant entire function and $a, b\in\mathbb{C}$ such that $b\not=a$. If $f$ and $f^{(1)}$ share $a$ and $b$ CM, then $f\equiv f^{(1)}$.
\end{theo2A}

In 2007, Bergweiler and Langley \cite{1a} established a relationship between $f^{(1)}(z)$ and $\Delta f(z)=f(z+1)-f(z)$ specifically showing that $f^{(1)}(z)\sim \Delta f(z)$, $r\rightarrow \infty$, $r\not\in E_1$, where $f(z)$ is a transcendental meromorphic function with $\rho(f)<1$ and $E_1$ is an $\varepsilon$-set.

\smallskip
A natural question arises: What is the uniqueness result when $f^{(1)}$ shares values with $\Delta_c f$ for a finite-order meromorphic function $f$? 

In response to this question, Qi et al. \cite{10} obtained the following results.

\begin{theo2B}\cite[Theorem 1.2]{10} Let $f$ be a meromorphic function of finite order. Suppose $f^{(1)}$ and $\Delta f$ share $a_1,a_2,a_3,a_4$ IM, where $a_1,a_2,a_3,a_4$ are four distinct finite values. Then $f^{(1)}\equiv \Delta_c f$.
\end{theo2B}

\begin{theo2C}\cite[Theorem 1.3]{10} Let $f$ be a transcendental meromorphic function such that its order of growth is not an integer or infinite and $a_1$ and $a_2$ be two distinct finite values. If $f^{(1)}$ and $\Delta_cf$  share $a_1, a_2, \infty $ CM, then $f^{(1)}\equiv \Delta_cf$.
\end{theo2C}

For further investigation, Qi et al. \cite{10} proposed the following question in the same paper:\par

{\bf Question 2.A.}  If we remove the condition that ``{\it order of growth is not an integer or infinite}'', does Theorem 2.C remain valid? \par

\smallskip
In the paper, we note that the difference analogue of the logarithmic derivative lemma is crucial to the demonstration of Theorem 2.C. However, for the difference logarithmic derivative lemma for infinite order meromorphic functions, a limited work has been done (see \cite{FLSY}, \cite{LFY1}), as far as the authors are aware. Recently, the last two authors of the present paper \cite{11c} proved that Theorem 2.C holds for a finite-order meromorphic function and obtained the following result.
\begin{theo2D}\cite[Theorem 2.1]{11c} Let $f$ be a non-constant meromorphic function such that $\rho(f)<\infty$ and let $a_1$ and $a_2$ be two distinct finite values. If $f^{(1)}$ and $\Delta_cf$  share $a_1, a_2, \infty $ CM, then $f^{(1)}\equiv \Delta_cf$.
\end{theo2D}

Therefore, the problem of investigating Theorem 2.D for an infinite-order meromorphic function appears to be worth considering. In this paper, we not only address this problem and provide a result that solves Question 2.A, but also improves Theorem 2.D. We now state our result:
\begin{theo}\label{t1} Let $f$ be a non-constant meromorphic function such that $\rho_1(f)<\infty$ and let $a_1$ and $a_2$ are two distinct finite values. If $f^{(1)}$ and $\Delta_cf$ share $a_1$, $a_2$ and $\infty$ CM, then $f^{(1)}\equiv \Delta_cf$.
\end{theo}

For entire function, we have the following result.

\begin{cor}\label{c1} Let $f$ be a non-constant entire function such that $\rho_1(f)<\infty$ and let $a_1$ and $a_2$ be two distinct finite values. If $f^{(1)}$ and $\Delta_cf$ share $a_1$ and $a_2$ CM, then $f^{(1)}\equiv \Delta_cf$.
\end{cor}

\begin{rem} The number of shared values cannot be reduced to two in Theorem \ref{t1}. For example, $f(z)=e^z$ and $e^c=3$. Note that $\Delta_cf(z)=2e^z$. Clearly $f^{(1)}$ and $\Delta_cf$ share $0$ and $\infty$ CM, but $f^{(1)}\not\equiv \Delta_cf$.
\end{rem}

Obviously, the scope of Theorem 2.D is extended in Theorem \ref{t1} by considering a larger class of meromorphic function $f$ such that
$\rho_1(f)=\infty$ and so Theorem \ref{t1} is an improvement of Theorem 2.D. But, our method does not work if $f$ is a non-constant meromorphic function such that $\rho_1(f)=\infty$. Therefore we ask the following open question: 

\begin{ques} If $f$ is a non-constant meromorphic function such that $\rho_1(f)=\infty$, is it then possible to establish Theorem \ref{t1}?\end{ques}

\medskip
We need the following lemmas for the proof of Theorem \ref{t1}.

\begin{lem}\label{l1} \cite[Theorem 3]{2} If $f$ and $g$ are two non-constant rational functions that share two values CM and one value IM, then $f\equiv g$.
\end{lem}

\begin{lem}\label{l2}\cite{8a} Let $f$ be a non-constant meromorphic function and $R(f)=\frac{P(f)}{Q(f)}$ where $P(f)=\sum_{k=0}^{p} a_{k}f^{k}$ and $Q(f)=\sum_{j=0}^{q} b_{j}f^{j}$ are two mutually prime polynomial in $f$. If $a_{k}, b_{j}\in S(f)$ are such that $a_{p}\not\equiv 0$ and $b_{q}\not\equiv 0$, then 
\[T(r,R(f))=\max\{p,q\}\;T(r,f)+S(r,f).\] 
\end{lem}

\begin{lem}\label{l3} \cite[Theorem 1.64]{14} Let $f_1,\ldots, f_n$ be non-constant meromorphic functions and $f_{n+1},\ldots, f_{n+m}$ be meromorphic functions such that $f_k\not\equiv 0\;\;(k=n+1,\ldots, n+m)$ and $\sum_{i=1}^{n+m} f_i\equiv A\in\mathbb{C}\setminus\{0\}$. If there exists a subset $I\subseteq \mathbb{R}^+$ satisfying $\text{mes}\;I=+\infty$ such that
\[\sum\limits_{i=1}^{n+m} N(r,0;f_i)+(n+m-1)\sum\limits_{\substack{i=1, i\neq j}}^{n+m} \ol N(r,f_i)<(\lambda+o(1))T(r,f_j)\;\;(r\to \infty, r\in I, j=1,2,\ldots,n),\] 
where $\lambda<1$, then there exist $t_i\in\{0,1\}\;(i=1,2,\ldots,m)$ such that $\sum_{i=1}^mt_if_{n+i}\equiv A$.
\end{lem}

\begin{lem}\label{l4} \cite[Lemma 6]{5} Let $f_0,\ldots, f_p$ be non-constant meromorphic functions such that $\Theta(\infty,f_j)=1$, $j = 0, \ldots, p$ and $\sum_{j=0}^p\lambda_jf_j\equiv 1$, where $\lambda_j\in\mathbb{C}\setminus\{0\}$. Then
\[\sideset{}{_{j=0}^p}{\sum}\Theta (0,f_j)\leq p+1-\frac{1}{p}.\]
\end{lem}

\begin{lem}\label{l6}\cite[Lemma 2.4]{9a} Let $f$ be a non-constant meromorphic function such that $\rho_1(f)<1$ and $c \in \mathbb{C}\setminus \{0\}$. Then 
\[T(r,f(z+c))=T(r,f)+S(r,f)\;\;\text{and}\;\;N(r,f(z+c))=N(r,f)+S(r,f).\]
\end{lem}

\begin{lem}\label{l5}\cite[Corollary 1]{13} Let $f$ be a non-constant meromorphic function on $\mathbb{C}$, and let $a_l,\ldots,a_q$ be distinct meromorphic functions on $\mathbb{C}$. Assume that $a_i$ are small functions with respect to $f$ for all $i=1,\ldots,q$. Then we have the second main theorem,
\[(q-2-\varepsilon)\;T(r, f) \leq \sideset{}{_{i=1}^{q}}{\sum}\ol N(r,a_i;f)+\varepsilon T(r, f),\]
for all $\varepsilon>0$ and for all $r\not\in E\subset (0,+\infty)$ such that $\int_E d \log \log r <+\infty.$
\end{lem}

\medskip
\begin{proof}[{\bf Proof of Theorem \ref{t1}}] We know that $f^{(1)}$ and $\Delta_cf$ share $a_1, a_2,\infty$ CM.
Now by the fundamental theorems, we get
\bea\label{nd1} T(r,f^{(1)})&\leq& \ol N(r,f^{(1)})+\ol N(r,a_1;f^{(1)})+\ol N(r,a_2;f^{(1)})+S(r,f^{(1)})\\&\leq&
\ol N(r,\Delta_cf)+\ol N(r,a_1;\Delta_cf)+\ol N(r,a_2;\Delta_cf)+S(r,f^{(1)})\nonumber\\&\leq& 3T(r,\Delta_cf)+S(r,f^{(1)}).\nonumber\eea

Similarly we have 
\bea\label{nd2} &&T(r,\Delta_cf)\leq 3T(r,f^{(1)})+S(r,\Delta_cf).\eea

Clearly from (\ref{nd1}) and (\ref{nd2}), we get $S(r,\Delta_cf)=S(r,f^{(1)})$.

\medskip
First we suppose $f$ is a non-constant rational function. Then $f^{(1)}$ is also a non-constant rational function.
Now from (\ref{nd2}), we see that $\Delta_cf$ is also a non-constant rational function. Then by Lemma \ref{l1}, we have $f^{(1)}\equiv \Delta_cf$.

\medskip
Next we suppose that $f$ is a transcendental meromorphic function. Then $f^{(1)}$ is also a transcendental meromorphic function. Now from (\ref{nd1}), we see that $\Delta_cf$ is also a transcendental meromorphic function.
Since $f^{(1)}$ and $\Delta_cf$ share $a_1, a_2, \infty$ CM, there exist two entire functions $P$ and $Q$ such that 
\bea\label{nd3} \frac{\Delta_cf-a_1}{f^{(1)}-a_1}=e^{P}\eea
and
\bea\label{nd4}\frac{\Delta_cf-a_2}{f^{(1)}-a_2}=e^{Q}.\eea

If $e^{P}\equiv e^{Q}$, then (\ref{nd3}) and (\ref{nd4}) give $\frac{\Delta_cf-a_1}{f^{(1)}-a_1}\equiv \frac{\Delta_cf-a_2}{f^{(1)}-a_2}$ and so $f^{(1)}\equiv \Delta_cf$. Henceforth we suppose $e^{P}\not\equiv e^{Q}$.
Now from (\ref{nd2}) and (\ref{nd3}), we get
\[T(r,e^{P})\leq T(r,\Delta_cf)+T(r,f^{(1)})+O(1)\leq 4T(r,f^{(1)})+S(r,f^{(1)})\leq 8T(r,f)+S(r,f)\]
and so $\rho(P)=\rho_1(e^P)\leq \rho_1(f)$. Since $\rho_1(f)<+\infty$, we have $\rho(P)=\rho_1(e^P)<+\infty$. Similarly we have $\rho(Q)=\rho_1(e^Q)<+\infty$.
We consider following three cases.\par

\medskip
{\bf Case 1.} Let $P$ be a polynomial and $Q$ be a transcendental entire function. If $e^P\equiv 1$ then from (\ref{nd3}), we get $f^{(1)}\equiv \Delta_cf$. Hence we suppose $e^P\not\equiv 1$.
Note that $\rho(e^P)=\deg(P)<+\infty$ and $\rho(e^Q)=+\infty$. Consequently $T(r,e^P)=S(r,e^Q)$. Note that 
\[T(r,Q^{(1)})=m(r,Q^{(1)})=m\left(r,(e^{Q})^{(1)}/e^{Q}\right)=S(r,e^{Q}).\]

Now solving (\ref{nd3}) and (\ref{nd4}) for $f^{(1)}$, we get 
\beas \label{k8} f^{(1)}=\frac{a_1e^{P}-a_2e^{Q}+a_2-a_1}{e^{P}-e^{Q}}\eeas
and so
\bea \label{k9} f^{(1)}(z+c)=\frac{a_1e^{P(z+c)}-a_2e^{Q(z+c)}+a_2-a_1}{e^{P(z+c)}-e^{Q(z+c)}}.\nonumber\eea

Therefore 
\bea\label{nd5} \Delta_cf^{(1)}=(a_2-a_1)\frac{e^{P}-e^{P+Q(z+c)}-e^{Q}+e^{P(z+c)+Q}-e^{P(z+c)}+e^{Q(z+c)}}{(e^{P}-e^{Q})(e^{P(z+c)}-e^{Q(z+c)})}.\eea 

Again from (\ref{nd3}) and (\ref{nd4}), we get
\bea\label{nd6}\Delta_cf=\frac{a_2e^{P}-a_1e^{Q}+(a_1-a_2)e^{P+Q}}{e^{P}-e^{Q}}.\eea

Differentiating (\ref{nd6}), we have
\bea\label{nd7} \Delta_cf^{(1)}=\frac{(a_1-a_2)\left((P^{(1)}-Q^{(1)})e^{P+Q}+Q^{(1)}e^{2P+Q}-P^{(1)}e^{P+2Q}\right)}{\left(e^{P}-e^{Q}\right)^2}.\eea

\smallskip
Let $g=e^P$. Then $g\not\equiv 1$ and $T(r,g)=S(r,e^Q)$.

\smallskip
Denote by $\ol N(r,g;e^Q\mid g\not=1)$ the reduced counting function of those zeros of $e^Q-g$ which are not the zeros of $g-1$. Also we denote by $\ol N(r,g;e^Q\mid g=1)$ the reduced counting function of those zeros of $e^Q-g$ which are the zeros of $g-1$.

\medskip
Suppose $\infty$ is not a Picard exceptional value of $f$.

Let $z_p$ be a pole of $f^{(1)}$ of multiplicity $p$ such that $g(z_p)\neq 1$. Then from (\ref{k8}), we see that $z_p$ is a zero of $e^{Q}-g$ of multiplicity $p$. Since $f^{(1)}$ has no simple poles, it follows that $p\geq 2$. Consequently 
\[\ol N(r,g;e^Q\mid g\neq 1)\leq \frac{1}{2} N(r,g;e^Q\mid g\neq 1)\leq \frac{1}{2} N(r,g;e^Q)\leq \frac{1}{2}T\left(r,e^Q\right)+S\left(r,e^Q\right).\]

On the other hand, we have 
\[\ol N(r,g;e^Q\mid g=1)\leq N(r,1;g)\leq T(r,g)=S(r,e^Q).\]

Therefore
\[\ol N(r,g;e^Q)=\ol N(r,g;e^Q\mid g\neq 1)+\ol N(r,g;e^Q\mid g=1)\leq \frac{1}{2}T\left(r,e^Q\right)+S\left(r,e^Q\right).\]

Now by Lemma \ref{l5}, we have 
\beas T(r,e^Q)\leq \ol N(r,g;e^{Q})+S(r,e^Q)\leq \frac{1}{2}T(r,e^Q)+S(r,e^Q),\eeas
which is impossible.\par

\medskip
Suppose $\infty$ is a Picard exceptional value of $f$.

Let $z_0$ be the zero of $e^Q-g$, i.e., $e^{Q(z_0)}=g(z_0)$. Now from (\ref{k8}), one can easily conclude that $z_0$ is a zero of $a_1g-a_2e^{Q}+a_2-a_1$, which shows that $g(z_0)=1$. Consequently 
\[\ol N(r,g;e^Q)\leq N(r,1;g)\leq S(r,e^Q).\]

Now by Lemma \ref{l5}, we get $T(r,e^Q)=S(r,e^Q)$, which is impossible. 

\medskip
{\bf Case 2.} Let $P$ is a transcendental entire function and $Q$ be a polynomial. In this case also we get $f^{(1)}\equiv \Delta_cf$, if we proceed in the same way as done in Case 1.\par

\medskip
{\bf Case 3.} Let $P$ and $Q$ be either both polynomials or both transcendental entire functions. 

Let $P$ and $Q$ be both constants. Set $e^{P}=c_1$ and $e^{Q}=d_1$. Then from (\ref{nd3}) and (\ref{nd4}), we have 
$\frac{\Delta_cf-a_1}{f^{(1)}-a_1}=c_1$ and $\frac{\Delta_cf-a_2}{f^{(1)}-a_2}=d_1$.
If $c_1=d_1$, then $f^{(1)}\equiv \Delta_cf$. Next let $c_1\neq d_1$. Clearly $f^{(1)}=\frac{a_1c_1-d_1a_2-a_1+a_2}{c_1-d_1}$ and so $f$ is a polynomial, which is absurd.\par

\smallskip
Henceforth we assume that $P$ and $Q$ are either both non-constant polynomials or both transcendental entire functions.\par

\medskip
$(A)$ Suppose first $P$ and $Q$ are both transcendental entire functions.

Now from (\ref{nd5}) and (\ref{nd7}), we get
\beas && e^{P+Q}-e^{2Q}-e^{P+Q+Q(z+c)}+e^{P(z+c)+2Q(z)}-e^{P(z+c)+Q}+e^{Q(z+c)+Q}-e^{2P}+e^{P+Q}\\
&&+e^{2P+Q(z+c)}-e^{P(z+c)+P+Q}+e^{P+P(z+c)}-e^{Q(z+c)+P}\\
&=&(P^{(1)}-Q^{(1)})e^{P(z+c)+P+Q}+Q^{(1)}e^{2P+Q+P(z+c)}-P^{(1)}e^{P+2Q+P(z+c)}\\
&&-(P^{(1)}-Q^{(1)})e^{Q(z+c)+P+Q}-Q^{(1)}e^{2P+Q+Q(z+c)}+P^{(1)}e^{P+2Q+Q(z+c)}
\eeas
and so
\bea\label{nd8} &&2e^{-R}-e^{-2R}+b_2e^{Q(z+c)-R}+e^{P(z+c)-2R}-e^{P(z+c)-R-P}+e^{Q(z+c)-R-P}+e^{Q(z+c)}\nonumber\\
&&+b_1e^{P(z+c)-R}+e^{P(z+c)-P}-e^{Q(z+c)-P}+b_3e^{P(z+c)+Q}+b_4e^{P(z+c)+Q-R}\nonumber\\
&&+b_5e^{Q(z+c)+Q}+b_6e^{Q(z+c)+Q-R}=1,\eea
where 
\bea\label{snd} 
b_1=-1-P^{(1)}+Q^{(1)}, b_2=-1+P^{(1)}-Q^{(1)}, b_3=-b_5=-Q^{(1)}, b_4=-b_6=P^{(1)}.\eea

We now prove that $R$ is not a polynomial. If not, suppose $R$ is a polynomial, say $p_1$. Since $e^P\not\equiv e^Q$, it follows that $e^{p_1}\not\equiv 1$ and so $e^{p_1(z+c)}\not\equiv 1$. Then (\ref{nd8}) gives
\bea\label{nd9} &&(b_2e^{-p_1}+1)e^{Q(z+c)}+(b_1+e^{-p_1})e^{-p_1}e^{P(z+c)}+(e^{-p_1}-1)(e^{-p_1(z+c)}-1)e^{P(z+c)-P(z)}\nonumber\\
&&+(b_5+b_6e^{-p_1})(e^{-p_1(z+c)}-1)e^{P(z+c)+Q(z)}=A_1,\eea
where $A_1=1-2e^{-p_1}+e^{-2p_1}$.\par

\smallskip
If $A_1\equiv 0$, then we get $(e^{-p_1}-1)^2=0$, i.e, $e^{p_1}=1$, i.e., $e^P\equiv e^Q$, which is a contradiction. Hence $A_1\not\equiv 0$. 
If possible suppose $b_2e^{-p_1}+1\equiv 0$. Now $b_2e^{-p_1}+1\equiv 0$ implies that 
$e^{p_1}\equiv -p_1^{(1)}+1$ and so by Lemma \ref{l2}, we get $T(r,e^{p_1})=S(r,e^{p_1})$, which shows that $p_1$ is a constant. Then we get $e^{p_1}=1$, which is impossible. Hence $b_2e^{-p_1}+1\not\equiv 0$. Similarly we can prove that $b_1+e^{-p_1}\not\equiv 0$.

\medskip
Fist suppose that $b_5+b_6e^{-p_1}\neq 0$.\par

\smallskip
Set $g_1=\frac{b_2e^{-p_1}+1}{A_1} e^{Q(z+c)}$, $g_2=\frac{(b_1+e^{-p_1})e^{-p_1}}{A_1} e^{P(z+c)}$,
$g_3=\frac{(e^{-p_1}-1)(e^{-p_1(z+c)}-1)}{A_1}e^{P(z+c)-P(z)}$ and $g_4=\frac{(b_5+b_6e^{-p_1})(e^{-p_1(z+c)}-1)}{A_1}e^{P(z+c)+Q(z)}$. Note that
\[\sideset{}{_{i=1}^{4}}{\sum} N(r,0;g_i)+3\sideset{}{_{\substack{i=1,i\neq j}}^{4}}{\sum} \ol N(r,g_i)<(\lambda+o(1))T(r,g_j)\;\;(r\to \infty, r\in I, j=1,2),\] 
$\lambda<1$. Then by Lemma \ref{l3}, there exist $t_i\in\{0,1\}, i=3,4$ such that $t_3g_3+t_4g_4=1$, i.e., 
\bea\label{nd10} t_3(e^{-p_1}-1)(e^{-p_1(z+c)}-1)e^{P(z+c)-P}+t_4(b_5+b_6e^{-p_1})(e^{-p_1(z+c)}-1)e^{P(z+c)+Q}=A_1.\eea

Clearly from (\ref{nd10}), we get $(t_3,t_4)\neq (0,0)$. 
Let $t_3=0$. Then $t_4=1$ and so from (\ref{nd9}) and (\ref{nd10}), we get
\bea\label{nd11} -(e^{-p_1}-1)(e^{-p_1(z+c)}-1)e^{-P}=(b_2e^{-p_1}+1)e^{-p_1(z+c)}+(b_1+e^{-p_1})e^{-p_1}.\eea

Now using Lemma \ref{l2} to (\ref{nd11}), we get $T(r,e^P)=S(r,e^P)$, which is impossible.

Let $t_4=0$. Then $t_3=1$ and so from (\ref{nd9}) and (\ref{nd10}), we get
\beas -(b_{5}+b_{6}e^{-p_{1}})(e^{-p(z+c)}-1)e^{Q}=(b_2e^{-p_1}+1)e^{-p_1(z+c)}+(b_1+e^{-p_1})e^{-p_1}\eeas
and so by Lemma \ref{l2}, we get a contradiction.

Hence $t_3=t_4=1$. Then from (\ref{nd9}) and (\ref{nd10}), we get
\beas (b_2e^{-p_1}+1)e^{Q(z+c)}+(b_1+e^{-p_1})e^{-p_1}e^{P(z+c)}=0\eeas
and so 
\beas e^{-p_1(z+c)}=-\frac{(b_1+e^{-p_1})e^{-p_1}}{b_2e^{-p_1}+1}.\eeas

Now using Lemmas \ref{l2} and \ref{l6}, we get $2T(r,e^{p_1})=T(r,e^{p_1})+S(r,e^{p_1})$, which shows that $p_1$ is a constant. Let $e^{-p_1}=d$. Then $d\neq 1$ and $P^{(1)}\equiv Q^{(1)}$.
Now applying Lemma \ref{l4} to (\ref{nd10}), we conclude that at least one of $P(z+c)-P$ and $Q+P(z+c)$ is a polynomial. Let $P(z+c)-P=p_2$. Then from (\ref{nd10}), we get
\bea\label{nd11a} dP^{(1)}e^{p_2}e^{2P}=e^{p_2}-1.\eea

Note that $T(r,P^{(1)})=S(r,e^P)$ and $T(r,e^{p_2})=S(r,e^P)$. Now using Lemma \ref{l2} to (\ref{nd11a}), we get $2T(r,e^P)=S(r,e^P)$, which is impossible.
Similarly if $Q+P(z+c)$ is a polynomial then we get a contradiction.\par

\medskip
Next suppose that $b_5+b_6e^{-p_1}\equiv 0$. Then from (\ref{nd9}), we get
\bea\label{nd11b} &&(b_2e^{-p_1}+1)e^{Q(z+c)}+(b_1+e^{-p_1})e^{-p_1}e^{P(z+c)}+(e^{-p_1}-1)(e^{-p_1(z+c)}-1)e^{P(z+c)-P(z)}\nonumber\\
&&=A_1.\eea

Now using Lemma \ref{l4} to (\ref{nd11b}), we conclude that $P(z+c)-P$ is a polynomial. Let $P(z+c)-P=p_4$. Then from (\ref{nd11b}), we get
\bea\label{nd12} &&(b_2e^{-p_1}+1)e^{Q(z+c)}+(b_1+e^{-p_1})e^{-p_1}e^{P(z+c)}=A_2,\eea
where $A_2=A_1-(e^{-p_1}-1)(e^{-p_1(z+c)}-1)e^{p_4}.$
If $A_2\not\equiv 0$, then using Lemma \ref{l4} to (\ref{nd12}), we get a contradiction. Hence $A_2\equiv 0$. Then from (\ref{nd12}), we get
\bea\label{nd12a} e^{p_1(z+c)}=\frac{(p_1^{(1)}-1)e^{-p_1}+1}{-e^{-2p_1}+(p_1^{(1)}+1)e^{-p_1}}.\eea

Now using Lemmas \ref{l2} and \ref{l6} to (\ref{nd12a}), we get $2T(r,e^{p_1})=T(r,e^{p_1})+S(r,e^{p_1})$, which shows that $p_1$ is a constant. Clearly $e^{p_1}\neq 1$ and $P^{(1)}\equiv Q^{(1)}\not\equiv 0$. Now $b_5+b_6e^{-p_1}\equiv 0$ implies that $e^{p_1}=1$, which is impossible. \par

\medskip
Hence $R$ is not a polynomial. We see from (\ref{snd}) that $b_1\not\equiv 0$ and $b_2\not\equiv 0$.
We can write (\ref{nd8}) in the following form 
\bea\label{n0} &&e^{-P}-e^{-Q}+b_5e^{P}+b_6e^{Q}+2e^{-Q(z+c)}+e^{R}+b_1e^{R(z+c)}+e^{-R+R(z+c)}-e^{R(z+c)-P}\nonumber\\
&&+e^{R(z+c)-Q}-e^{R-Q(z+c)}-e^{-R-Q(z+c)}+b_3e^{R(z+c)+P}+b_4e^{Q+R(z+c)}=-b_2.\eea

We consider following sub-cases.\par

\medskip
{\bf Sub-case 3.1.} Let $-R+R(z+c)=q_1$, where $q_1$ is a polynomial. Then (\ref{n0}) gives 
\bea\label{n1} &&e^{-P}-(1+e^{q_1})e^{-Q}+(b_5+b_4e^{q_1})e^{P}+b_6e^{Q}+2e^{-Q(z+c)}+(1+b_1e^{q_1})e^{R}\nonumber\\
&&+e^{q_1}e^{P-2Q}-e^{-q_1}e^{P(z+c)-2Q(z+c)}-e^{q_1}e^{-P(z+c)}+b_3e^{q_1}e^{2P-Q}=-b_2-e^{q_1}.\eea

Now we divide two possible sub-cases:\par

\medskip
{\bf Sub-case 3.1.1.} Let $-b_2-e^{q_1}\not\equiv 0$. We prove that $P-2Q$ and $2P-Q$ are not polynomials.\par

\smallskip 
Let $P-2Q=q_2$. Then from (\ref{n1}), we have 
\bea\label{n2}&&(e^{-q_2}-e^{-q_1-2q_2+q_2(z+c)})e^{-2Q}+(2e^{-q_1-q_2+q_2(z+c)}-(1+e^{q_1}))e^{-Q}\\
&&+(b_6+b_3e^{q_1+q_2}+(1+b_1e^{q_1}))e^Q+(b_5+b_4e^{q_1})e^{q_2}e^{2Q}=-b_2-e^{q_1}-e^{q_1+q_2}+e^{-q_1+q_2(z+c)}.\nonumber\eea

If atleast one of $e^{-q_2}-e^{-q_1-2q_2+q_2(z+c)}$, $2e^{-q_1-q_2+q_2(z+c)}-(1+e^{q_1})$, $b_5+b_3e^{q_1+q_2}+(1+b_1e^{q_1})$, $b_6+b_4e^{q_1}$, $-b_2-e^{q_1}-e^{q_1+q_2}+e^{q_1+q_2(z+c)}$ is not identically zero, then using Lemma \ref{l2} to (\ref{n2}), we get $T(r,e^Q)=S(r,e^Q)$, which is a contradiction. Hence 
\beas &&e^{-q_2}-e^{-q_1-2q_2+q_2(z+c)}\equiv 0,\;\;2e^{-q_1-q_2+q_2(z+c)}-1-e^{q_1}\equiv 0,\\&&b_6+b_3e^{q_1+q_2}+1+b_1e^{q_1}\equiv 0,\;\;b_5+b_4e^{q_1}\equiv 0,\\&&-b_2-e^{q_1}-e^{q_1+q_2}+e^{-q_1+q_2(z+c)}\equiv 0.\eeas 

Now from $e^{-q_2}-e^{-q_1-2q_2+q_2(z+c)}\equiv 0$ and  $2e^{-q_1-q_2+q_2(z+c)}-(1+e^{q_1})\equiv 0$ we get $e^{q_1}\equiv -1$. Then $b_6+b_4e^{q_1}\equiv 0$ gives $Q^{(1)}-P^{(1)}\equiv 0$, i.e., $Q-P\equiv c\in\mathbb{C}$. Since $P-2Q=q_2$, we deduce that $-Q=c+q_2$, which is impossible.\par

\smallskip
Let $Q-2P=q_3$. Then from (\ref{n1}), we get
\beas && ((1+b_1e^{q_1})e^{-p_3}-e^{q_1-3q_3+q_3(z+c)}-e^{q_1-q_3+q_3(z+c)}+1)e^{-P}\\
&&+((1+e^{q_1})e^{-q_3}+2e^{-2q_3+q_3(z+c)})e^{-2P}+e^{q_1-2q_3}e^{-3P}+(b_5+b_4e^{q_1})e^P+b_6e^{q_3}e^{2P}\\
&&=-b_2-e^{q_1}+b_3e^{q_1-q_3}\eeas
and by Lemma \ref{l2}, we get $T(r,e^P)=S(r,e^P)$, which is a contradiction.\par

\medskip
{\bf Sub-case 3.1.2.} Let $-b_2-e^{q_1}\equiv 0$. Then from (\ref{snd}), we get $1-R^{(1)}=e^{-R+R(z+c)}$. This shows that $1$ is a Picard exceptional value of $R^{(1)}$ and so we let $R^{(1)}=1+e^{R_1}$, where $R_1$ is an entire function. Since $-R+R(z+c)=q_1$, we have $-R^{(1)}+R^{(1)}(z+c)=q_1^{(1)}$ and so $-e^{R_1}+e^{R_1(z+c)}=q_1^{(1)}$, i.e., $e^{R_1}+q^{(1)}_1=e^{R_1(z+c)}$. Let $q^{(1)}_1\not\equiv 0$. Then $N(r,-q_1^{(1)};e^{R_1})=0$. Now by Lemma \ref{l5}, we get $T(r,e^{R_1})=S(r,e^{R_1})$, which is absurd. Hence $q^{(1)}_1\equiv 0$ and so $q_1$ is a constant, say $k$. Then $1-R^{(1)}=e^k$ and so $R$ is a polynomial, which is absurd.\par

\medskip
{\bf Sub-case 3.2.} Let $R(z+c)-P=q_4$, where $q_4$ is a polynomial. Then (\ref{n0}) gives
\bea\label{n3}&& e^{-P}-e^{-Q}+(b_5+b_1e^{q_4})e^P+(b_6+e^{q_4})e^Q+2e^{-Q(z+c)}+(1+e^{q_4})e^R\nonumber\\
&& +b_3e^{q_4}e^{2P}+b_4e^{q_4}e^{P+Q}-e^{R-Q(z+c)}-e^{-R-Q(z+c)}=-b_2+e^{q_4}.\eea

Now we divide the following sub-cases.\par

\medskip
{\bf Sub-case 3.2.1.} Let $-b_2+e^{q_4}\not \equiv 0$. We prove that $P+Q$ is not a polynomial.\par

\smallskip
Let $P+Q=q_5$, where $q_5$ is a polynomial. Now  from (\ref{n3}), we get
\beas && 2e^{-q_5}e^{-P(z+c)}+((b_6+e^{q_4})e^{2q_5+q_4}+e^{q_4+q_5})e^{-2P(z+c)}-e^{2q_4+3q_5}e^{-5P(z+c)}\\
&&((b_5+b_1e^{q_4})e^{-q_4-q_5}-e^{-q_4-2q_5})e^{2P(z+c)}-e^{-2q_4-3q_5}e^{3P(z+c)}\\
&&+((1+e^{q_4})e^{-2q_4-3q_5}+b_3e^{-q_4-2q_5})e^{4P(z+c)}=-b_2+e^{q_4}-b_4e^{q_4+q_5}\eeas
and so by Lemma \ref{l2}, we conclude that $T(r,e^{P(z+c)})=S(r,e^{P(z+c)})$, which is a contradiction. Hence $P+Q$ is not a polynomial.\par

\smallskip
Now by Lemma \ref{l3}, there exists $s_i\in\{0,1\}$, $(i=1,2)$ such that 
\bea\label{n4} && -s_1e^{R-Q(z+c)}-s_2e^{-R-Q(z+c)}=-b_2+e^{q_4}.\eea
Let $s_2=0$. Then from (\ref{n4}), we get $-s_1e^{R-Q(z+c)}=-b_2+e^{q_4}$, which shows that $R-Q(z+c)$ is a polynomial, say $q_6$. Now from (\ref{n3}), we get 
\bea\label{n5}&& -e^R-(b_5+b_1e^{q_4})e^{2P}+(b_6+e^{q_4})e^{P+Q}+2e^{q_6}e^Q+(1+e^{q_4})e^{2P-Q}\nonumber\\
&& +b_3e^{q_4}e^{3P}+b_4e^{q_4}e^{2P+Q}-e^{q_6}e^{P-2Q}=-1.\eea

If possible suppose $2P-Q$ is a polynomial, say $q_7$. We know that $R(z+c)-P=q_4$ and $R-Q(z+c)=q_6$. Then we conclude that $P(z+c)=-q_6+q_7+q_4$, which gives a contradiction. Hence $2P-Q$ is not a polynomial. Similarly we can prove that $2P+Q$ and $P-2Q$ are not polynomials. Now using Lemma \ref{l4} to (\ref{n5}), we get a contradiction.\par

\smallskip
If $s_1=0$, then proceeding in the same way as done above, we get a contradiction. Hence $s_1=s_2=1$. Then from (\ref{n4}), we get 
\[e^{-R-Q(z+c)}(e^{2R}+1)=-(-b_2+e^{q_4}).\]

From (\ref{snd}), we see that $T(r,-b_2+e^{q_4})=S(r,e^R)$. Then from above, we get $N(r,-1, e^{2R})=S(r,e^{2R})$ and so by the second fundamental theorem, we get a contradiction.\par

\medskip
{\bf Sub-case 3.2.2.} Let $-b_2+e^{q_4}\equiv 0$. Then from (\ref{n3}), we get 
\bea\label{n6}&& e^{-2P-Q}-e^{-P-2Q}+(b_5+b_1e^{q_4})e^{-Q}+(b_6+e^{q_4})e^{-P}+2e^{-Q(z+c)-P-Q}\nonumber\\
&&+(1+e^{q_4})e^{-2Q}+b_3e^{q_4}e^R-e^{-2Q-Q(z+c)}-e^{-2P-Q(z+c)}=-b_4e^{q_4}.\eea

Now we prove that $2P+Q$, $P+2Q$ and $P+Q+Q(z+c)$ are not polynomials.\par

\smallskip
Let $2P+Q=q_8$, where $q_8$ is a polynomial. Then from (\ref{n6}), we get 
\beas && (b_3e^{q_8-2q_4-3q_8(z+c)}-e^{2q_8-3q_4-3q_8(z+c)})e^{9P(z+c)}+(b_6+b_1e^{q_4})e^{q_8-2q_4-2q_8(z+c)}e^{6P(z+c)}\\
&&+e^{2q_4+q_8(z+c)}e^{-3P(z+c)}+2e^{-q_4+q_8-q_8(z+c)}e^{5P(z+c)}+(1+e^{q_4})e^{2q_8-4q_4-4q_8(z+c)}e^{12P(z+c)}\\
&&-e^{2q_8-4q_4-5q_8(z+c)}e^{14P(z+c)}-e^{2q_4+q_8(z+c)}e^{-5P(z+c)}=-e^{-q_8}-b_4e^{q_4}\eeas
and by Lemma \ref{l2}, we get $T(r,e^{P(z+c)})=S(r,e^{P(z+c)})$, which is a contradiction.\par

\smallskip
Hence $2P+Q$ is not a polynomial.

If $P+2Q$ is a polynomial, then proceeding in the same way as done above, we get a contradiction. Hence $P+2Q$ is not a polynomial.\par

\smallskip
Let $P+Q+Q(z+c)=q_9$, where $q_9$ is polynomial. Now from (\ref{n6}), we get 
\bea\label{n7}&& e^{-2P-Q}-e^{-P-2Q}+(b_5+b_1e^{q_4})e^{-Q}+(b_6+e^{q_4})e^{-P}\nonumber\\
&&+(1+e^{q_4})e^{-2Q}+(b_3e^{q_4}-e^{q_9})e^R-e^{-q_9}e^{-R}=b_7,\eea
where $b_7=-2e^{-q_9}-b_4e^{q_4}.$ If $b_7\not\equiv 0$, then using Lemma \ref{l3} to (\ref{n7}), we get a contradiction. Hence $b_7\equiv 0$. Then from (\ref{n7}), we get
\bea\label{n8}&& e^{-P-2Q}-e^{-3Q}+(b_5+b_1e^{q_4})e^{P-2Q}+(b_6+e^{q_4})e^{-Q}\nonumber\\
&&+(1+e^{q_4})e^{P-3Q}+(b_3e^{q_4}-e^{q_9})e^{2R}=e^{-q_9}.\eea

If possible suppose $P-2Q=q_{10}$, where $q_{10}$ is a polynomial. We know that $R(z+c)-P=q_4$ and $P+Q+Q(z+c)=q_9$. Then by a simple calculation, we get $-3Q=-q_4+q_9-2q_{10}+q_{10}(z+c)$, which is absurd. Hence $P-2Q$ is not a polynomial. Similarly we can prove that $P-3Q$ is not a polynomial. Now using Lemma \ref{l4} to (\ref{n8}), we get a contradiction.

Hence $2P+Q$, $P+2Q$ and $P+Q+Q(z+c)$ are not polynomials.
Now by Lemma \ref{l3}, there exists $t_i\in\{0,1\},(i=1,2)$ such that 
\bea\label{n9} -t_1e^{-2Q-Q(z+c)}-t_2e^{-2P-Q(z+c)}=-b_4e^{q_4}\not\equiv 0.\eea

Let $t_2=0$. Then from (\ref{n9}), we get $-e^{-2Q-Q(z+c)}=-b_4e^{q_4}$, which shows that $-2Q-Q(z+c)$ is a polynomial, say $q_{11}$. Now from (\ref{n6}) and (\ref{n9}), we get 
\bea\label{n10}&& e^{-P-2Q}-e^{-3Q}+(b_5+b_1e^{q_4})e^{P-2Q}+(b_6+e^{q_4})e^{-Q}\nonumber\\
&&+(1+e^{q_4})e^{P-3Q}+b_3e^{q_4}e^{2R}-e^{q_{11}}e^{-R}=-2e^{q_{11}}.\eea

It is easy to verify that $P+2Q$, $P-2Q$ and $P-3Q$ are transcendental entire functions. Then using Lemma \ref{l4} to (\ref{n10}), we get a contradiction.\par

\smallskip
Let $t_1=0$. Then proceeding in the same way as done above, we get a contradiction. Hence $t_1=t_2=1$. Now from (\ref{n9}), we have 
\[e^{-2P-Q(z+c)}(e^{2R}+1)=b_4e^{q_4}.\]

From (\ref{snd}), we see that $T(r,b_4)=S(r,e^R)$. Then from above we get $N(r,-1,e^{2R})=S(r,e^{2R})$ and so by the second fundamental theorem, we get a contradiction.\par

\medskip
{\bf Sub-case 3.3.} Let $R(z+c)+P$ be a polynomial. We omit the proof since the same can be carried out in the line of proof of Sub-case 3.2.\par

\medskip
{\bf{Sub-case 3.4.}} Let $R(z+c)-Q=q_{12}$, where $q_{12}$ is a polynomial. Then (\ref{n0}) gives
\bea\label{0n11} &&e^{-P}-e^{-Q}+b_5e^{P}+b_6e^{Q}+2e^{-Q(z+c)}+e^{R}+b_1e^{q_{12}}e^{Q}+e^{q_{12}}e^{-P+2Q}-e^{q_{12}}e^{-R}\nonumber\\
&&-e^{R-Q(z+c)}-e^{-R-Q(z+c)}+b_3e^{q_{12}}e^{P+Q}+b_4e^{q_{12}}e^{2Q}=-b_2-e^{q_{12}}.\eea

Now we divide the following sub-cases.\par

\medskip
{\bf{Sub-case 3.4.1.}} Let $-b_2-e^{q_{12}}\not\equiv 0$. First we prove that $P+Q$ and $-P+2Q$ are transcendental. \par

\smallskip
Let $P+Q=q_{13}$, where $q_{13}$ is a polynomial. Since $R(z+c)-Q=q_{12}$, we get $2Q(z+c)-P=(q_{13}(z+c)-q_{13}-q_{12})/2$. Then from (\ref{0n11}), we get 
\beas &&(1+b_6e^{q_{13}}+b_1e^{q_{13}+q_{12}})e^{-P}+(b_5-e^{-q_{13}})e^P+e^{-q_{13}}e^{2P}+(b_4-1)e^{q_{13}+q_{12}}e^{-2P}\\&&+(2-e^{-q_{13}}e^{2P}-e^{q_{13}+q_{12}}e^{-2P})e^{(q_{13}(z+c)-q_{13}-q_{12})/2}e^{P/2}+e^{q_{12}+q_{13}}e^{-3P}\\
&&=-b_2-e^{q_{12}}-b_3e^{q_{12}+q_{13}}.\eeas

Using Lemma \ref{l2}, we get a contradiction. Hence $P+Q$ is not a polynomial.\par

\smallskip
Let $-P+2Q=q_{14}$, where $q_{14}$ is a polynomial. Then from (\ref{0n11}), we get 
\beas &&e^{q_{14}}e^{-2Q}-(1+e^{q_{14}+q_{12}})e^{-Q}+(b_6+2e^{-q_{14}(z+c)+q_{12}}+e^{-q_{14}}+b_1e^{q_{12}})e^Q\\&&+(b_5e^{-q_{14}}-e^{q_{14}+q_{12}-q_{14}(z+c)}+b_4e^{q_{12}})q^{2Q}+b_3e^{q_{12}-q_{14}}e^{3Q}\\&&=-b_2-e^{q_{12}}-e^{q_{14}+q_{12}}+e^{q_{14}+q_{12}+q_{14}(z+c)}\eeas

Using Lemma \ref{l2}, we get a contradiction. Hence $-P+2Q$ is not a polynomial. Similarly we can show that $P+2Q$ is not a polynomial.\par

\smallskip
Hence $P+Q$ and $-P+2Q$ are transcendental. Now we apply Lemma \ref{l3} to (\ref{0n11}). Then there exist $u_i\in\{0,1\}\;(i=1,2)$ such that
\bea\label{n12}-u_1e^{R-Q(z+c)}-u_2e^{-R-Q(z+c)}=-b_2-e^{q_{12}}.\eea

Clearly $(u_1,u_2)\neq(0,0)$. Let $u_2=0$. Then from (\ref{n12}), we get $-e^{R-Q(z+c)}=-b_2-e^{q_{12}}$ and so $R-Q(z+c)$ is a polynomial, say $q_{15}$. Therefore (\ref{n11}) gives
\bea\label{n13}&&e^{-P-2Q}-e^{-3Q}+b_5e^{P-2Q}+(b_6+b_1e^{q_{12}})e^{-Q}+(2e^{q_{15}}-e^{q_{12}})e^{-P-Q}\\\nonumber&&+e^{P-3Q}+e^{q_{12}}e^{-P}-e^{q_{15}}e^{-2P}+b_3e^{q_{12}}e^{P-Q}=-b_4e^{q_{12}}.\eea

Using Lemma \ref{l4} to (\ref{n13}), we get a contradiction. Hence $u_2\neq 0$.\par

\smallskip
Let $u_1=0$. Then from (\ref{n12}), we get $-e^{-R-Q(z+c)}=-b_2-e^{q_{12}}$ and so $-R-Q(z+c)$ is a polynomial. Therefore proceeding in the same way as done above, we get a contradiction.\par

\smallskip
Hence $u_1=u_2=1$. Then from (\ref{n12}), we get $N(r,1,e^{-2R})=S(r,e^{-2R})$ and so using second fundamental theorem, we get $T(r,e^{-2R})=S(r,e^{-2R})$, which is a contradiction.\par

\medskip
{\bf{Sub-case 3.4.2.}}  Let $-b_2-e^{q_{12}}\equiv 0$. Then from (\ref{0n11}), we have
\bea\label{n13} &&e^{-Q}-e^{P-2Q}+b_5e^{2P-Q}+b_6e^{P}+2e^{R-Q(z+c)}+e^{2R}+b_1e^{q_{12}}e^{P}\nonumber\\
&&+e^{q_{12}}e^{Q}-e^{2R-Q(z+c)}-e^{-Q(z+c)}+b_3e^{q_{12}}e^{2P}+b_4e^{q_{12}}e^{P+Q}=e^{q_{12}}.\eea

Let $P+Q=q_{16}$, where $q_{16}$ is a polynomial. Then from (\ref{n13}), we have 
\beas &&(b_6e^{q_{16}}+b_1e^{q_{16}+q_{12}}+1)e^{-Q}+(b_5e^{2q_{16}}-1)e^{-3Q}+\\
&&(2e^{q_{16}}e^{-2Q}-e^{2q_{16}}e^{-4Q}-1)e^{(q_{16}(z+c)-q_{12})/2}e^{Q/2}\\&&+e^{2q_{16}}e^{-4Q}+e^{q_{12}}e^Q+b_3e^{2q_{16}+q_{12}}e^{-2Q}=e^{q_{12}}-b_4e^{q_{16}+q_{12}}.\eeas

Using Lemma \ref{l2}, we get a contradiction. Hence $P+Q$ is not a polynomial.\par

\smallskip
Let $P-2Q=q_{17}$, where $q_{17}$ is a polynomial. Then from (\ref{n13}), we get 
\beas &&(1-e^{q_{17}(z+c)-q_{12}})e^{-Q}+(b_5e^{2q_{17}}+b_4e^{q_{17}+q_{12}})e^{3Q}+(b_6e^{q_{17}}+b_1e^{q_{17}+q_{12}}+e^{2q_{17}})e^{2Q}\\&&+(e^{q_{12}}-e^{q_{17}(z+c)-q_{12}+2q_{17}})e^Q+b_3e^{2q_{17}+q_{12}}e^{4Q}=e^{q_{12}}+e^{q_{17}}.\eeas

Using Lemma \ref{l2}, we get a contradiction. Hence $P-2Q$ is not a polynomial. 

Similarly we can show that $2P-Q$ is not a polynomial.\par

\smallskip
Hence $P+Q$, $P-2Q$ and $2P-Q$ are transcendental.
Now we Lemma \ref{l3} to (\ref{n13}). Then there exist $v_i\in\{0,1\}\;(i=1,2)$ such that 
\bea\label{n14} v_12e^{R-Q(z+c)}-v_2e^{2R-Q(z+c)}=e^{q_{12}}.\eea

Clearly $(v_1,v_2)\neq(0,0)$. Let $v_2=0$. Then (\ref{n14}) gives $2e^{R-Q(z+c)}=e^{q_{12}}$ and so $R-Q(z+c)$ is a polynomial. Now proceeding in the same way as done in sub-case 3.4.1, we get a contradiction. Hence $v_2\neq 0$.\par

\smallskip
Let $v_1=0$. Then (\ref{n14}) gives $-e^{2R-Q(z+c)}=e^{q_{12}}$ and so $2R-Q(z+c)$ is a polynomial. Now proceeding in the same way as done in sub-case 3.4.1, we get a contradiction.\par

\smallskip
Hence $v_1=v_2=1$. Then (\ref{n14}) gives $N(r,-2,e^{R})=S(r,e^{R})$ and so using second fundamental theorem, we get $T(r,e^{R})=S(r,e^{R})$, which is a contradiction.\par

\medskip
{\bf Sub-case 3.5.} Let $R(z+c)+Q$ be a polynomial. We omit the proof since the same can be carried out in the line of proof of Sub-case 3.4.\par

\medskip
Finally from {\bf Sub-case 3.1-Sub-case 3.5}, we prove that $-R+R(z+c)$, $R(z+c)-P$, $R(z+c)+P$, $R(z+c)+Q$ and $R(z+c)-Q$ are transcendental entire functions. Now applying Lemma \ref{l3} to (\ref{n0}). Then there exists $w_i\in\{0,1\}(i=1,2)$ such that 
\bea\label{n11}-w_1e^{R-Q(z+c)}-w_2e^{-R-Q(z+c)}=-b_2.\eea

We divide the following sub-cases.\par

\medskip
{\bf Sub-case 3.A.} Let $w_2=0$. Then from (\ref{n11}), we conclude that $Q(z+c)-R$ is a polynomial, say $q_{18}$. Now (\ref{n0}) gives that 
\bea\label{s1} &&e^{q_{18}}e^{-Q}-e^{q_{18}}e^{P-2Q}+b_5e^{q_{18}}e^{2P-Q}+b_6e^{q_{18}}e^P+e^{q_{18}}e^{2R}+b_1e^{P(z+c)}+e^{q_{18}}e^{R(z+c)}\nonumber\\
&&-e^{P(z+c)-P}+e^{P(z+c)-Q}-e^{R}+b_3e^{P(z+c)+P}+b_4e^{-q_{18}} e^{-P+2Q}=-2.\eea 

Let $P-2Q=q_{19}$, where $q_{19}$ is a polynomial. Then from (\ref{s1}), we get 
\beas && e^{q_{18}}e^{-Q}+(e^{q_{19}(z+c)+2q_{18}+2q_{19}}+e^{2q_{18}+q_{19}+q_{19}(z+c)}-e^{q_{19}})e^Q\\
&&+(b_6+e^{q_{19}})e^{q_{18}+q_{19}}+b_1e^{q_{19}(z+c)+2q_{19}+q_{18}})e^{2Q}+b_5e^{q_{18}+2q_{19}}e^{3Q}+b_3e^{q_{19}(z+c)+2q_{18}+3q_{19}}e^{4Q}\\
&&=e^{2q_{18}+q_{19}+q_{19}(z+c)}-b_4e^{-q_{18}-q_{19}}+e^{q_{18}+q_{19}}-2\eeas
and by Lemma \ref{l2}, we get $T(r,e^Q)=S(r,e^Q)$, which is a contradiction.\par

\smallskip
Hence $P-2Q$ is not a polynomial.

If $2P-Q$ is a polynomial, then  proceeding in the same way as done above, we get a contradiction. Hence $2P-Q$ is not a polynomial.\par

\smallskip
Let $P(z+c)-Q=q_{20}$ where $q_{20}$ is a polynomial. Then from (\ref{s1}), we get
\bea \label{s2}&&b_1e^{q_{20}}e^{Q}+e^{q_{18}}e^{-Q}+(e^{q_{20}}-e^{q_{18}})e^{P-2Q}+b_5e^{q_{18}}e^{2P-Q}+b_6e^{q_{18}}e^P+e^{q_{18}}e^{2R}\nonumber\\
&&-e^{q_{20}}e^{-R}-e^{-q_{18}}e^{-3R}+b_3e^{q_{20}}e^{P+Q}+b_4e^{-q_{18}} e^{-P+2Q}=-2-e^{q_{20}}.\eea

If possible suppose $P+Q$ is a polynomial, say $q_{21}$. We know that $Q(z+c)-R=q_{18}$ and $-P+P(z+c)=q_{20}$. Then by a simple calculation, we get$-P=q_{18}+q_{20}-q_{21}-q_{21}(z+c)$, which gives a contradiction. Hence $P+Q$ is not a polynomial. If $2+e^{q_{20}}\not\equiv 0$, then using Lemma \ref{l4} to (\ref{s2}), we get a contradiction. Hence $2+e^{q_{20}}\equiv 0$. Now from (\ref{s2}), we get 
\bea \label{s3}&&b_1e^{q_{20}}e^{-R}+e^{q_{18}}e^{-P-Q}+(e^{q_{20}}-e^{q_{18}})e^{-2Q}+b_5e^{q_{18}}e^{R}+e^{q_{18}}e^{P-2Q}\nonumber\\
&&-e^{q_{20}}e^{-2P+Q}-e^{-q_{18}}e^{-4P+3Q}+b_3e^{q_{20}}e^{Q}+b_4e^{-q_{18}} e^{-2R}=-b_6e^{q_{18}}.\eea

It is easy to verify that $-4P+3Q$ and  $-P-Q$ are transcendentals. Since $b_6\not \equiv 0$, using Lemma \ref{l4} to (\ref{s3}), we get a contradiction.\par

\smallskip
Therefore $P(z+c)-Q$ is not a polynomial.

Hence we prove that $P-2Q$, $2P-Q$ and $P(z+c)-Q$ are transcendental entire functions. Now we apply Lemma \ref{l3} to (\ref{s1}). Then there exists $v_i\in\{0,1\},(i=1,2)$ such that 
\bea\label{s4} -v_1e^{P(z+c)-P}+v_2b_3e^{P(z+c)+P}=-2.\eea
 
Let $v_2=0$. Then from (\ref{s4}), we get $e^{P(z+c)-P}=2$ which shows that $-P+P(z+c)=q_{22}$, where $q_{22}$ is a polynomial. Then from (\ref{s1}), we get
\bea\label{s5} &&e^{q_{22}}e^{-R}+e^{q_{18}}e^{-P-Q}-e^{q_{18}}e^{-2Q}+b_5e^{q_{18}}e^{R}+e^{q_{18}}e^{P-2Q}\nonumber\\
&&+e^{q_{22}}e^{-Q}-e^{-q_{18}}e^{-4P+3Q}+b_3e^{q_{22}}e^{P}+b_4e^{-q_{18}} e^{-2R}=-b_6e^{q_{18}}.\eea

If possible suppose $-4P+3Q$ be a polynomial, say $q_{23}$. We know that $Q(z+c)-R=q_{18}$ and $-P+P(z+c)=q_{22}$. Then by a simple calculation, we get $-5P=-3q_{18}+4q_{22}+q_{23}+q_{23}(z+c)$, which gives a contradiction. Hence $-4P+3Q$ is not polynomial. Similarly $P+Q$ is not polynomial. Now using Lemma \ref{l4} to (\ref{s5}), we get a contradiction. \par

\smallskip
Let $v_1=0$. Then  proceeding in the same way as above, we get a contradiction. Hence $v_1=v_2=1$. Now from (\ref{s4}), we get 
$e^{P(z+c)-P}(b_3e^{2P}-1)=-2$,
which shows that $N(r,1;b_3e^{2P})=S(r,b_3e^{2P})$ and so by Lemma \ref{l5}, we get a contradiction. \par

\medskip
{\bf Sub-case 3.B.} Let $w_1=0$. Then proceeding in the same way as done in sub-case 3.A, we get a contradiction.\par

\medskip
{\bf Sub-case 3.C.} Let $w_1=w_2=1$. Then from (\ref{n11}), we get 
$e^{-R-Q(z+c)}(e^{2R}+1)=b_2$, 
which shows that $N(r,-1,e^{2R})=S(r,e^{2R})$ and so by Lemma \ref{l5}, we get a contradiction.\par

\medskip
$(B)$ Suppose $P$ and $Q$ are both non-constant polynomials.

In this case we use the same methodology as used in the proof of Case 3 but with a little modification as follows:

We apply Lemma \ref{l4} for $p_i$'s as non-constants instead of $p_i$'s as non-polynomials. 

This completes the proof. 
\end{proof}

\section{\bf{Periodicity on a transcendental entire function}}

We know that, if $f(z)$ is a periodic function with period $c$, then $f^{(k)}(z)$ is also a periodic function with the same period $c$, but the converse is not true in general. As for example, let $f(z)=\sin z+z$. Clearly $f^{(1)}(z)$ is a periodic function with period $2\pi$, but $f(z)$ is not periodic with period $2\pi$. However if $f^{(k)}(z)$ is a periodic function of period $c$, then $f(z)$ can be expressed as $f(z)=\varphi(z)+p(z)$, where $\varphi(z)$ is a periodic function of period $c$ and $p(z)$ is a polynomial such that $\deg(p)\leq k$. Related to the periodicity of entire functions, Wang and Hu \cite{11} mentioned the following conjecture posed by Yang.

\smallskip
{\bf Yang's Conjecture.} Let $f$ be a transcendental entire function and $k$ be a positive integer. If $ff^{(k)}$
is a periodic function, then $f$ is also a periodic function.\par

\smallskip
In 2018, Wang and Hu \cite{11a} proved that {\bf Yang's Conjecture} is true for $k=1$ by giving the following result.
\begin{theo3A}\cite[Theorem 1.1]{11a} Let $f$ be a transcendental entire function and $k$ be a positive integer. If $(f^2)^{(k)}$
is a periodic function, then $f$ is also a periodic function.
\end{theo3A}

In 2019, Liu and Yu \cite{8} further improved Theorem 3.A in a direction where $f^2$ is replaced by a more general polynomial  $a_nf^n+a_{n-1}f^{n-1}+\ldots+a_1f$, $a_1, \ldots, a_n(\neq 0)$ are constants and $n\geq 2$. Now we recall their result.

\begin{theo3B}\cite[Theorem 1.7]{8} Let $f$ be a transcendental entire function such that $\rho_1(f)<1$ and $N(r,0;f)=S(r,f)$ and let $k$ be a positive integer and $n\geq 2$. If $(a_1f+\ldots+a_nf^n)^{(k)}$ is a periodic function, then $f$ is also a
periodic function, where $a_n\neq 0$.
\end{theo3B}

Recently Wei et al. \cite{11b} proved that Theorem 3.B holds without the conditions ``$\rho_1(f)<1$ and $N(r,0;f)=S(r,f)$''.
Also in the same paper, Wei et al. \cite{11b} considered the problem of periodicity of $f^n(f^m-1)f^{(k)}$. Actually they obtained the following results.

\begin{theo3C}\cite[Corollary 1.5]{11b} Assume that $f$ is a transcendental entire function. If
$f^n(f^m-1)f^{(1)}$ is a periodic function, then $f$ is also a periodic function.
\end{theo3C}

\begin{theo3D}\cite[Theorem 1.6]{11b} Let $f$ be a transcendental entire function and $n,\;m$ and
$k$ be positive integers. If $f^n(f^m-1)f^{(k)}$ is a periodic function with period $c$,
and one of following conditions is satisfied:
\begin{enumerate}
\item[(i)] $f(z)=e^{h(z)}$, where $h(z)$ is an entire function;
\item[(ii)] $f(z)$ has a non-zero Picard exceptional value and $f(z)$ is of finite
order;
\item[(iii)] $f^n(f^m-1)f^{(k+1)}$ is a periodic function with period $c$,
\end{enumerate}
then $f(z)$ is also a periodic function.
\end{theo3D}

Motivated by {\bf Yang's Conjecture} and Theorem 3.D, Wei et al. \cite{11b} raised the following question for further research.\par

\smallskip
{\bf Question 3.A.} Let $f$ be a transcendental entire function and $m,\;n$ and 
$k$ be positive integers. If $\frac{f^{(k)}}{f^n(f^m-1)}$ is a periodic function with period $c$, is it true that $f$ is
also a periodic function?\par

\smallskip
In the same paper, Wei et al. \cite{11b} proved that Question 3.A is solvable for the case $m=n=k=1$. 

In the paper, inspired by Theorem 3.D we will try to give an affirmative answer of Question 3.A. Now we state our result.
\begin{theo}\label{t3.1} Let $f$ be a transcendental entire function and $n,\;m$ and
$k$ be positive integers. If $\frac{f^{(k)}}{f^n(f^m-1)}$ is a periodic function with period $c$,
and one of following conditions is satisfied:
\begin{enumerate}
\item[(I)] $f(z)=e^{P(z)}$, where $P(z)$ is a non-constant entire function;
\item[(II)] $f(z)$ has a non-zero Picard exceptional value and $\rho_1(f)<+\infty$;
\item[(III)] $\frac{f^{(k+1)}}{f^n(f^m-1)}$ is a periodic function with period $c$,
\end{enumerate}
then $f(z)$ is also a periodic function.
\end{theo}

Following  lemma is required in the proof of Theorem \ref{t3.1}.
\begin{lem}\label{l3.1}\cite[Theorem 2]{ss1} If $P(z)$ is a non-constant polynomial and $f(z)$ is an entire
function which is not periodic, then $P(f(z))$ cannot be periodic either.
\end{lem}

\smallskip
\begin{proof}[{\bf Proof of Theorem \ref{t3.1}}] 
We know $\frac{f^{(k)}}{f^n(f^m-1)}$ is a periodic function with period $c$ and so 
\bea\label{k1} \frac{f^{(k)}(z+c)}{f^n(z+c)(f^m(z+c)-1)}=\frac{f^{(k)}(z)}{f^n(z)(f^m(z)-1)}.\eea

\smallskip
{\bf (I).} Let $f(z)=e^{P(z)}$, where $P(z)$ is a non-constant entire function. Differentiating $k$-times, we get $f^{(k)}(z)=\tilde P (P)e^{P(z)}$, where $\tilde P(P)$ is a differential polynomial in $P^{(1)}$. Since $f$ is a transcendental entire function, we get $\tilde P(P)\not\equiv 0$. Note that $P^{(i)}\in S(f)$ for $i=1,2,\ldots$. Therefore $\tilde P(P)\in S(f)$ and $\tilde P(P(z+c))\in S(f(z+c))$. Now from (\ref{k1}), we get
\bea\label{k1a} \frac{\tilde P(P(z+c))}{e^{(n-1)P(z+c)}(e^{mP(z+c)}-1)}=\frac{\tilde P(P(z))}{e^{(n-1)P(z)}(e^{mP(z)}-1)}.\eea 

Then using Lemma \ref{l2} to (\ref{k1a}), we get 
\[T(r,e^{P(z)})+S(r,e^{P(z)})=T(r,e^{P(z+c)})+S(r,e^{P(z+c)})\]
and so $S(r,e^{P(z)})=S(r,e^{P(z+c)})$. Let 
\[\alpha(z)=\frac{\tilde P (P(z))}{\tilde P (P(z+c))}.\]

Obviously $\alpha\not\equiv 0$. Then (\ref{k1a}) gives
\bea\label{k2} \alpha(z) e^{(n+m-1)P(z+c)-(n-1)P(z)}-\alpha(z) e^{(n-1)(P(z+c)-P(z))}-e^{mP(z)}=-1.\eea

Set 
\[g_1(z)=\alpha e^{(n+m-1)P(z+c)-(n-1)P(z)},\;\; g_2(z)=\alpha e^{(n-1)(P(z+c)-P(z))}\;\; \text{and}\;\; g_3(z)=e^{mP(z)}.\]

Clearly $T(r,\alpha)=S(r,e^P)=S(r,g_3)$. Note that 
\beas \sideset{}{_{i=1}^3}{\sum} N(r,0,g_i)+2\sideset{}{_{i=1}^2}{\sum}\ol N(r,g_i)<(\lambda+o(1))T(r,g_3)\;\;(r\to \infty, r\in I,\; \lambda<1).\eeas

Now we apply Lemma \ref{l3} to (\ref{k2}). Then there exists $t_i\in \{0,1\}(i=1,2)$ such that  
\bea\label{k3} t_1\alpha(z) e^{(n+m-1)P(z+c)-(n-1)P(z)}-t_2\alpha(z) e^{(n-1)(P(z+c)-P(z))}=-1.\eea

Let $t_2=0$. Now from (\ref{k3}), we get 
\[\alpha(z)e^{(n+m-1)P(z+c)-(n-1)P(z)}=-1\]
and so from (\ref{k2}), we obtain 
\[\alpha(z)e^{(n-1)P(z+c)-(n+m-1)P(z)}=-1.\]

Then $e^{mP(z+c)}=e^{mP(z)}$, i.e., $f^m(z+c)=f^m(z)$ and so $f(z+c)=tf(z)$, $t^m=1$. It means that $f$ is a periodic function with period $mc$.\par

\smallskip
Let $t_1=0$. Then proceeding in the same way as above, we get $f(z)$ is a periodic function with period $mc$. 

Let $t_1=t_2=1$. Then from (\ref{k3}), we get 
\[\alpha(z)e^{(n-1)(P(z+c)-P(z))}(e^{mP(z+c)}-1)=-1,\]
which shows that 
\[N(r,1,e^{P(z+c)})\leq N(r,\alpha)\leq T(r,\alpha)= S(r,e^{mP(z+c)}).\]

Now by Lemma \ref{l5}, we get $T(r,e^{mP(z+c)})=S(r,e^{mP(z+c)})$, which is a contradiction.\par

\smallskip
{\bf (II).} Let $a$ be a non-zero Picard exceptional value of $f$ and $\rho_1(f)<+\infty$. Then there exists an entire function $h(z)$ such that 
\bea\label{k4} f(z)=a+e^{h(z)}.\eea

Now differentiating (\ref{k4}) $k$-times, we get 
\bea\label{k5} f^{(k)}(z)=\hat P(h^{(1)}(z))e^{h(z)},\eea
where $\hat P(h^{(1)})=(h^{(1)})^k+\hat P_{k-1}(h^{(1)})$,
$\hat P_{k-1}$ is a differential polynomial in $h^{(1)}$ of degree less then or equal to $k-1$. Clearly $\hat P(h^{(1)})\not\equiv 0$. 
Now from (\ref{k1}), (\ref{k4}) and (\ref{k5}), we get 
\bea\label{k6} &&\frac{1}{\hat P(h^{(1)}(z+c))}(e^{(m+n)h(z+c)}+b_{m+n}e^{(m+n-1)h(z+c)}+..+b_{n+1}e^{nh(z+c)}+..+b_1)e^{h(z)}\nonumber\\
&&=\frac{1}{\hat P(h^{(1)}(z))}(e^{(m+n)h(z)}+b_{m+n}e^{(m+n-1)h(z)}+..+b_{n+1}e^{nh(z)}+..+b_1)e^{h(z+c)},\eea
where 
\[b_{m+n}=\binom{m+n}{1}a,\cdots,b_{n+1}=\binom{m+n}{m}a^{m}-1,\cdots, b_2=a^{n-1}((m+n)a^m-n), b_1=a^{m+n}-a^n.\]

Clearly $(b_1, b_{n+1})\neq (0,0)$. Now using Lemma \ref{l2} to (\ref{k6}), we get $S(r,e^{h(z+c)})=S(r,e^{h(z)})$. So 
\[T(r,\hat P(h^{(1)}(z)))=S(r,e^{h(z)})=S(r,e^{h(z+c)})\]
and 
\[T(r,\hat P(h^{(1)}(z+c)))=S(r,e^{h(z+c)})=S(r,e^{h(z)}).\]

Let 
\[G(z)=\frac{\hat P(h^{(1)}(z))}{\hat P(h^{(1)}(z+c))}.\]

Clearly $G\not\equiv 0$. Now (\ref{k6}) gives 
\bea\label{k7} && Ge^{(m+n-1)h(z+c)+h(z)}+b_{m+n}Ge^{(m+n-2)h(z+c)+h(z)}+\cdots+b_3e^{h(z+c)+h(z)}\nonumber\\&&+b_1Ge^{h(z)-h(z+c)}-e^{(m+n)h(z)}-b_{m+n}e^{(m+n-1)h(z)}-\cdots+(G-1)b_2e^{h(z)}=b_1.\eea

We now divide the following two cases.\par

\smallskip
{\bf Case 1.} Let $h$ be transcendental. Since $\rho_1(f)<+\infty$, we get $\rho(h)=\rho_1(e^h)<+\infty$.\par

\smallskip
First suppose $b_1=a^{m+m}-a^n\neq 0$. Using Lemma \ref{l4} to (\ref{k7}), we get at least one of $h(z)+jh(z+c)(j=1,2,\ldots,m+n-1)$ and $h(z)-h(z+c)$ is a polynomial. For the sake of simplicity, we assume $(m+n-1)h(z+c)+h(z)$ is a polynomial, say $p$. Then (\ref{k7}) gives
\bea\label{k8} && Ge^{p}+b_{m+n}Ge^{p}e^{-h(z+c)}+\cdots+b_3Ge^pe^{-(m+n-2)h(z+c)}+b_1Ge^{p}e^{-(m+n)h(z+c)}\nonumber\\
&&-e^{(m+n)p}e^{-(m+n)(m+n-1)h(z+c)}-b_{m+n}e^{(m+n-1)p}e^{-(m+n-1)(m+n-1)h(z+c)}-\cdots\nonumber\\
&&+(G-1)b_2e^pe^{-(m+n-1)h(z+c)}=b_1.\eea

Now using Lemma \ref{l2} to (\ref{k8}), we get $T(r,e^{h(z+c)})=S(r,e^{h(z+c)})$, which is impossible.\par 

\medskip
Next we suppose $b_1=a^{m+m}-a^n=0$. Then $a^m=1$. Since $b_2=a^{n-1}((m+n)a^m-n)$, it follows that $b_2\neq 0$. Now from (\ref{k7}), we get
\bea\label{k9} && Ge^{(m+n-1)h(z+c)}+b_{m+n}Ge^{(m+n-2)h(z+c)}+\cdots+b_3Ge^{h(z+c)}\nonumber\\
&&-e^{(m+n-1)h(z)}-b_{m+n}e^{(m+n-2)h(z)}-\cdots-b_3e^{h(z)}=-(G-1)b_2.\eea

Let $G\equiv 1$. Then from (\ref{k9}), we get 
\bea\label{k10} &&e^{(m+n-1)h(z+c)}+b_{m+n}e^{(m+n-2)h(z+c)}+\cdots+b_3e^{h(z+c)}\nonumber\\
&&=e^{(m+n-1)h(z)}+b_{m+n}e^{(m+n-2)h(z)}+\cdots+b_3e^{h(z)}. \eea

Set $f_1(z)=e^{h(z)}$ and $P(z)=z^{m+n-1}+b_{m+n}z^{m+m-2}+\cdots+b_3z$. Clearly $P(z)$ is non-constant. Also from (\ref{k10}), we see that $P(f_1(z))$ is periodic. Now using Lemma \ref{l3.1} to (\ref{k10}), we conclude that $f_1(z)$ must be a periodic function. Hence $f(z)$ is a periodic function.\par

\smallskip
Let $G\not \equiv 1$. Now using Lemma \ref{l3} to (\ref{k9}), we get a contradiction.\par

\smallskip
{\bf Case 2.} Let $h$ be a polynomial. In this case we use the same methodology as used in the proof of Case 1, but with a little modification as follows:

We apply Lemma \ref{l4} for $p$'s as non-constants instead of $p$'s as non-polynomials.\par

\smallskip
{\bf (III).} Let $\frac{f^{(k+1)}}{f^n(f^m-1)}$ be a periodic function with period $c$. Then by Theorem 1.9 \cite{11b}, we deduce that $f$ is also a periodic function.
This completes the proof.
\end{proof}

{\bf Compliance of Ethical Standards:}\par

{\bf Conflict of Interest.} The authors declare that there is no conflict of interest regarding the publication of this paper.\par

{\bf Data availability statement.} Data sharing not applicable to this article as no data sets were generated or analysed during the current study.

\end{document}